\documentclass[10pt]{article}
\def\R{\mathbb{R}}
\usepackage[margin=0.85in]{geometry}
\usepackage{amsmath,amssymb,amsthm,booktabs,hyperref,pgfplots,xcolor,caption,float}
\newtheorem{theorem}{Theorem}
\pgfplotsset{compat=1.18}
\usetikzlibrary{arrows.meta}
\colorlet{g}{green!45!black}
\colorlet{b}{blue!75!black}
\colorlet{k}{black!60}
\hypersetup{colorlinks=true,urlcolor=blue,linkcolor=blue,citecolor=blue}
\newcommand{\rev}[1]{{\color{black}#1}}
\title{Heterogeneous Hegselmann-Krause Dynamics Need Not Converge}
\author{S.\ Rasoul Etesami \qquad Behrouz Touri\\[4pt]
\small Department of Industrial and Enterprise Systems Engineering\\
\small Coordinated Science Lab\\
\small University of Illinois Urbana-Champaign\\
\small\texttt{\{etesami1,touri1\}@illinois.edu}}
\date{}
\begin{document}
\maketitle

\begin{abstract}
We give a six-agent counterexample to convergence of scalar heterogeneous
Hegselmann--Krause \rev{dynamics.}
\end{abstract}

\section{Introduction}
In bounded-confidence opinion dynamics, agents synchronously average the opinions of neighbors whose opinions lie within prescribed confidence ranges\rev{,} with the Hegselmann--Krause (HK) model as a prototypical example \cite{hk2002,krause2000}. In its homogeneous form all agents have the same confidence radius. Finite-time termination and quantitative convergence bounds for the homogeneous model are established in \cite{mohajer2013convergence,bhcn2013,eb2015}; the total-energy framework of \cite{chazelle2011} supplies related tools for bidirectional multiagent systems. For a survey of bounded-confidence models, see \cite{BERNARDO2024111302}.

In the heterogeneous HK model, agent $i$ has its own fixed confidence radius $\rho_i>0$ and observes agent $j$ when $|x_j-x_i|\le\rho_i$. Unlike the homogeneous model, this interaction relation need not be symmetric. Over the last decades there \rev{have} been many theoretical partial results to support the conjecture that the heterogeneous HK model always converges~\cite{mb2012,inertiaHK,suyu2016,eb2015}. Hegarty, Ognissanti and Wedin \cite{how2026} recently constructed counterexamples to two auxiliary conjectures and a theorem in \cite{mb2012}. This counterexample\rev{,} however, did not resolve the general convergence conjecture as the resulting trajectory \rev{converges}. Motivated \rev{by} and based on their work, here we give an explicit six-agent setting with a periodic solution for the heterogeneous scalar HK rule.

\section{The model, the nonexistence, and the existence of cycles}
Consider $n$ agents $[n]:=\{1,\ldots,n\}$ with initial opinions $x_i(0)\in\R$ and fixed observation radius
$\rho_i>0$ for each $i\in [n]$. The HK model prescribes  the evolution of each \rev{agent's} opinion $x_i(t)\in\mathbb R$ as a 
discrete-time  $t\in\mathbb Z_{\ge0}$ dynamics with the update rule
\begin{equation}\label{eq:HK}
  x_i(t+1)=\frac1{|N_i(t)|}\sum_{j\in N_i(t)}x_j(t),
\end{equation}
where $N_i(t)=\{j\in[n]:|x_j(t)-x_i(t)|\le\rho_i\}$.

Extensive numerical experiments and partial theoretical results have suggested that this dynamics always converges, i.e., for any $n$ and any vector of positive observation radii $\rho=(\rho_1,\ldots,\rho_n)^T$ and any initial \rev{condition} $x(0)\in \R^n$, the above dynamics converges. However, a formal proof of this conjecture has remained elusive. One possible reason is that nonconvergent behavior may emerge only in larger systems involving more agents and heterogeneous confidence bounds. In fact, in our initial attempt to establish convergence, we showed that the heterogeneous HK dynamics cannot enter a periodic orbit of length two, as formally stated in the following result.

\begin{theorem}\label{thm:no-2-cycle}
The scalar heterogeneous HK model has no nonconstant orbit of period two.
\end{theorem}
\begin{proof}
Suppose two distinct states $p,q$ satisfy $q=F(p)$ and $p=F(q)$, where
$F$ denotes the HK update. Choose $i$ with
$|q_i-p_i|=\Delta:=\max_j|q_j-p_j|>0$. Reflecting the line if necessary,
assume $q_i=p_i+\Delta$. Since $q_i$ averages opinions within $\rho_i$
of $p_i$, we have $\Delta\le\rho_i$.

Let $I=N_i(p)\cap N_i(q)$, $D=N_i(p)\setminus I$, and
$E=N_i(q)\setminus I$. A departing source cannot leave through the
right boundary, because
$q_j\le p_j+\Delta\le q_i+\rho_i$ for $j\in D$.
It therefore leaves through the left boundary, giving
$p_j\le q_j+\Delta<q_i-\rho_i+\Delta\le q_i$.
Similarly, a new source cannot enter from the left: if
$p_j<p_i-\rho_i$, then $q_j\le p_j+\Delta<q_i-\rho_i$.
Hence every $j\in E$ satisfies
$q_j\ge p_j-\Delta>p_i+\rho_i-\Delta\ge p_i$.

Removing the departing sources from the average for $q_i$ can only
increase it; removing the entering sources from the average for $p_i$
can only decrease it. Consequently,
\[
 \sum_{j\in I}p_j\ge |I|q_i,
 \qquad
 \sum_{j\in I}q_j\le |I|p_i.
\]
But $q_j\ge p_j-\Delta$ for all $j$, with strict inequality at $j=i$.
Since self inclusion gives $i\in I$,
\[
 \sum_{j\in I}q_j
 >\sum_{j\in I}p_j-|I|\Delta
 \ge |I|(q_i-\Delta)=|I|p_i,
\]
a contradiction.
\end{proof}

One may conjecture that such a periodic orbit does not exist for any orbit length $\ell$. However, this is not true\rev{,} as will be discussed shortly. 

Consider the solution $z^*$ of the fixed-point equation
\begin{equation}\label{eq:fixed}
 z^*=\frac{2161}{15552}+\frac{z^*}{1024}.
\end{equation}
This linear equation has a unique solution. Based on the solution to this equation, we have the following result.

\begin{table}[H]
\centering
\begin{tabular}{clll}
\toprule
Agent & Initial opinion & Radius & Role\\
\midrule
1 & $0$          & $1/1000$ & left endpoint\\
2 & $z^*=34576/248589$ & $1/4$ & left green, $g$\\
3 & $2/5$    & $59/100$ & left blue, $b$\\
4 & $3/5$    & $59/100$ & right blue, $1-b$\\
5 & $1-z^*=214013/248589$ & $1/4$ & right green, $1-g$\\
6 & $1$          & $1/1000$ & right endpoint\\
\bottomrule
\end{tabular}
\caption{Initial opinions and fixed confidence radii.}
\label{tab:config}
\end{table}

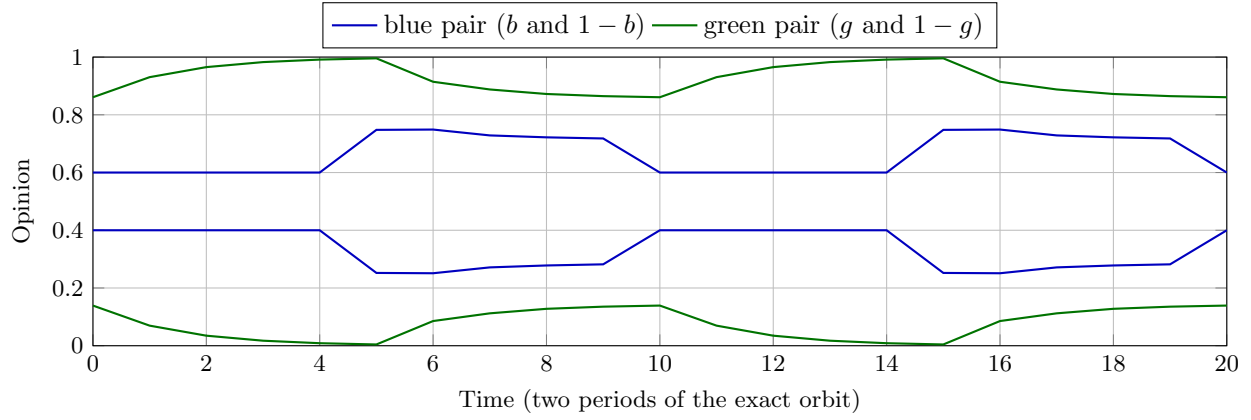
\begin{figure}[htbp]
\begin{center}
\begin{tikzpicture}
\begin{axis}[
 width=0.96\textwidth,height=5.4cm,
 xlabel={Time (two periods of the exact orbit)},
 ylabel={Opinion},xmin=0,xmax=20,ymin=0,ymax=1.0,
 xtick={0,2,4,6,8,10,12,14,16,18,20},ytick={0,0.2,0.4,0.6,0.8,1},grid=major,
 legend style={at={(0.5,1.03)},anchor=south,legend columns=2},
 tick label style={font=\small},label style={font=\small}
]
\addplot[b,thick] coordinates {(0,0.4) (1,0.4) (2,0.4) (3,0.4) (4,0.4) (5,0.252173265913) (6,0.251086632956) (7,0.271376649812) (8,0.278049436017) (9,0.281964532823) (10,0.4) (11,0.4) (12,0.4) (13,0.4) (14,0.4) (15,0.252173265913) (16,0.251086632956) (17,0.271376649812) (18,0.278049436017) (19,0.281964532823) (20,0.4)};
\addlegendentry{blue pair ($b$ and $1-b$)}
\addplot[b,thick,forget plot] coordinates {(0,0.6) (1,0.6) (2,0.6) (3,0.6) (4,0.6) (5,0.747826734087) (6,0.748913367044) (7,0.728623350188) (8,0.721950563983) (9,0.718035467177) (10,0.6) (11,0.6) (12,0.6) (13,0.6) (14,0.6) (15,0.747826734087) (16,0.748913367044) (17,0.728623350188) (18,0.721950563983) (19,0.718035467177) (20,0.6)};
\addplot[g,thick] coordinates {(0,0.13908901842) (1,0.06954450921) (2,0.034772254605) (3,0.0173861273025) (4,0.00869306365125) (5,0.00434653182562) (6,0.0855065992461) (7,0.112197744068) (8,0.127858131293) (9,0.135302522437) (10,0.13908901842) (11,0.06954450921) (12,0.034772254605) (13,0.0173861273025) (14,0.00869306365125) (15,0.00434653182562) (16,0.0855065992461) (17,0.112197744068) (18,0.127858131293) (19,0.135302522437) (20,0.13908901842)};
\addlegendentry{green pair ($g$ and $1-g$)}
\addplot[g,thick,forget plot] coordinates {(0,0.86091098158) (1,0.93045549079) (2,0.965227745395) (3,0.982613872698) (4,0.991306936349) (5,0.995653468174) (6,0.914493400754) (7,0.887802255932) (8,0.872141868707) (9,0.864697477563) (10,0.86091098158) (11,0.93045549079) (12,0.965227745395) (13,0.982613872698) (14,0.991306936349) (15,0.995653468174) (16,0.914493400754) (17,0.887802255932) (18,0.872141868707) (19,0.864697477563) (20,0.86091098158)};
\end{axis}
\end{tikzpicture}
\end{center}
\caption{Two periods of the exact orbit from the initial data in
Table~\ref{tab:config}, obtained by evaluating
Table~\ref{tab:itinerary} at $z=z^*$. Decimal plot coordinates are used
only for visualization. The endpoints stay at $0$ and $1$.}\label{fig:orbit}
\end{figure}

\begin{theorem}\label{thm:cycle}
Let $n=6$. For the initial condition $x(0)=(0,z^*,2/5,3/5,1-z^*,1)^T\in \R^n$ with observation radii profile $\rho=(0.001,0.25,0.59,0.59,0.25,0.001)$ (see Table~\ref{tab:config} and Fig.~\ref{fig:config}), the resulting heterogeneous HK dynamics \rev{is periodic with fundamental period ten}. As a result, the heterogeneous HK model does not always converge. 
\end{theorem}
\begin{proof}
For the given initial condition, we have two anchor agents $i=1$ and $i=6$ at $0$ and $1$, respectively (which we show to be fixed throughout the evolution). Call agents $2$ and $5$\rev{, initially located at $z^*$ and $1-z^*$,} the green agents, and the agents at $2/5$ and $1-2/5=3/5$ the blue agents. 

The symmetry of this initial condition around $1/2$ and the symmetry of the the observation radii implies symmetrical evolution of the dynamics around $1/2$ throughout. Assuming the endpoints remain fixed, write
$x(t)=\bigl(0,g(t),b(t),1-b(t),1-g(t),1\bigr)^\top$. We will prove simultaneously that this representation holds for all
time.

Consider the working intervals
\begin{equation}\label{eq:region}
 g\in I_g:=[\frac1{250},\frac7{50}],\qquad
 b\in I_b:=[\frac14,\frac25].
\end{equation}
Here $0<g<b<1-b<1-g<1$. Each endpoint observes only itself, since its
nearest neighbor is at distance $g\ge1/250>1/1000$. The left blue
agent always observes $0,g,b,1-b$, excludes $1$, and observes $1-g$ precisely
when $b+g\ge41/100$. The left green agent always observes $0,g$ and observes
$b$ precisely when $b-g\le1/4$. It excludes the other three agents:
\[
 1-b-g\ge\frac{46}{100}>\frac14,\qquad
 1-2g\ge\frac{72}{100}>\frac14,\qquad
 1-g>\frac14.
\]
The right agents have the corresponding reflected neighborhoods. For a generic $b\in I_b,g\in I_g$  define the \rev{binary indicator} variables
\begin{equation}\label{eq:gates}
 v=\mathbf1_{\{b+g\ge41/100\}},\qquad
 u=\mathbf1_{\{b-g\le1/4\}}.
\end{equation}
Note that $u,v$ depend on $b,g$. With this, $b(t)$ and $g(t)$ can be written as
\begin{equation}\label{eq:F}
 \begin{aligned}
 b(t+1)&=\begin{cases}
       2/5,&v(t)=1,\\
       (1+g(t))/4,&v(t)=0,
      \end{cases}
 &g(t+1)&=\begin{cases}
       (b(t)+g(t))/3,&u(t)=1,\\
       g(t)/2,&u(t)=0,
      \end{cases}
 \end{aligned}
\end{equation}
where $u(t),v(t)$ are the corresponding indicator variables to $b(t),g(t)$. 
For example, when $v(t)=1$ the \rev{opinions of the} five agents observed by the left blue agent
\rev{sum} to $0+g(t)+b(t)+(1-b(t))+(1-g(t))=2$, which explains the reset to $2/5$.

\begin{figure}[htbp]
\centering
\begin{tikzpicture}[
 x=8.6cm,y=0.76cm,
 grp/.style={circle,inner sep=0pt,draw=black!75,line width=0.4pt,fill opacity=0.72},
 stem/.style={black!45,densely dotted,line width=0.4pt},
 lab/.style={font=\scriptsize,align=center,inner sep=1.5pt},
 guide/.style={black!18,dashed,line width=0.35pt},
 bar/.style={line width=1.1pt,{Bar[width=5pt]}-{Bar[width=5pt]}},
 halo/.style={font=\scriptsize,fill=white,fill opacity=0.92,text opacity=1,inner sep=1.5pt},
 tick/.style={font=\tiny,black!65,fill=white,fill opacity=0.92,text opacity=1,inner sep=1pt}
]
\def\zg{0.139089018420}
\def\zgr{0.860910981580}
\foreach \x in {0,\zg,0.4,0.6,\zgr,1}
 \draw[guide] (\x,1.15) -- (\x,-6.10);
\draw[-{Latex[length=5pt]},line width=0.7pt] (-0.28,0) -- (1.28,0);
\node[font=\small] at (1.315,0) {$x$};
\foreach \x/\tk in {0/0,0.25/0.25,0.5/0.5,0.75/0.75,1/1}{
 \draw[black!65,line width=0.5pt] (\x,0.11) -- (\x,-0.11);
 \node[tick] at (\x,-0.33) {\tk};}
\foreach \x in {\zg,0.4,0.6,\zgr}
 \draw[stem] (\x,0.11) -- (\x,0.62);
\node[grp,fill=g,minimum size=6pt] at (\zg,0.75) {};
\node[grp,fill=b,minimum size=6pt] at (0.4,0.75) {};
\node[grp,fill=b,minimum size=6pt] at (0.6,0.75) {};
\node[grp,fill=g,minimum size=6pt] at (\zgr,0.75) {};
\node[lab,text=g] at (\zg,1.45) {$g$: $z^*$\\$\rho=1/4$};
\node[lab,text=b] at (0.4,1.45) {$b$: $0.4$\\$\rho=0.59$};
\node[lab,text=b] at (0.6,1.45) {$1-b$: $0.6$\\$\rho=0.59$};
\node[lab,text=g] at (\zgr,1.45) {$1-g$: $1-z^*$\\$\rho=1/4$};
\foreach \x in {0,1}{
 \draw[stem] (\x,-0.11) -- (\x,-0.83);
 \node[grp,fill=k,minimum size=6pt] at (\x,-0.95) {};}
\node[lab,text=k] at (0,-1.62) {$0$\\$\rho=1/1000$};
\node[lab,text=k] at (1,-1.62) {$1$\\$\rho=1/1000$};
\node[font=\scriptsize\itshape,text=k,anchor=west] at (0.30,-0.95)
 {fixed endpoint agents};
\draw[black!35,line width=0.4pt] (-0.28,-2.42) -- (1.28,-2.42);
\node[halo,text=b,anchor=west] at (-0.26,-3.03)
 {$b$: $[-0.19,\,0.99]$ includes five agents; the endpoint $1$ lies just outside};
\draw[bar,b] (-0.19,-3.35) -- (0.99,-3.35);
\fill[b] (0.4,-3.35) circle (1.6pt);
\node[halo,text=b,anchor=west] at (-0.26,-3.68)
 {$1-b$: interval $[0.01,\,1.19]$, the reflection};
\draw[bar,b] (0.01,-4.00) -- (1.19,-4.00);
\fill[b] (0.6,-4.00) circle (1.6pt);
\node[halo,text=g,anchor=west] at (-0.26,-4.33)
 {$g$: $[z^*-1/4,\,z^*+1/4]$ observes only $0$ and $g$};
\draw[bar,g] (-0.110910981580,-4.65) -- (0.389089018420,-4.65);
\fill[g] (\zg,-4.65) circle (1.6pt);
\node[halo,text=g,anchor=west] at (-0.26,-4.98)
 {$1-g$: interval $[3/4-z^*,\,5/4-z^*]$, the reflection};
\draw[bar,g] (0.610910981580,-5.30) -- (1.110910981580,-5.30);
\fill[g] (\zgr,-5.30) circle (1.6pt);
\node[halo,text=k,anchor=west] at (-0.26,-5.68)
 {endpoints: $\rho=1/1000$ (intervals too short to distinguish at this scale)};
\foreach \x in {0,1}
 \draw[k,line width=1.1pt] (\x,-5.92) -- (\x,-6.06);
\end{tikzpicture}
\caption{Initial opinions and confidence intervals. Each green agent
observes its nearest endpoint and itself; each blue agent observes five agents.
The configuration and the confidence radii are symmetric about $1/2$.}
\label{fig:config}
\end{figure}
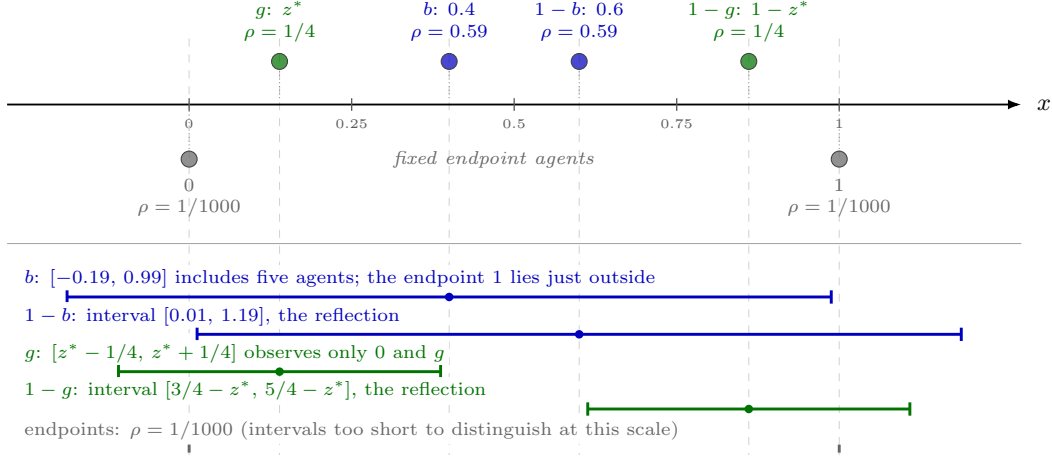

If the initial condition $z^*$ \rev{were} instead set to any $z\in I:=\left[\frac{139}{1000},\frac7{50}\right]$,
starting from $(b,g)=(2/5,z)$, successive applications of
\eqref{eq:F} give Table~\ref{tab:itinerary} due to the fact that 
\[
 \begin{aligned}
 b(t)+g(t)&>41/100 &&(t=0,1,2,3,9),\\
 b(t)+g(t)&<41/100 &&(t=4,5,6,7,8),\\
 b(t)-g(t)&>1/4   &&(t=0,1,2,3,4),\\
b(t)-g(t)&<1/4   &&(t=5,6,7,8,9).
 \end{aligned}
\]
These strict inequalities verify the switches throughout $I$.

Noticing that $z^*$, the \rev{fixed-point} solution of $z=2161/15552+z/1024$, belongs to the interior of $I$\rev{, which leads} to the corresponding \rev{period-$10$} orbit outlined in that table, concludes the proof. 

\begin{table}[H]
\centering
\begin{tabular}{clll}
\toprule
$t$ & $b(t)$ & $g(t)$ & $(v(t),u(t))$\\
\midrule
0 & $2/5$ & $z$ & $(1,0)$\\
1 & $2/5$ & $z/2$ & $(1,0)$\\
2 & $2/5$ & $z/4$ & $(1,0)$\\
3 & $2/5$ & $z/8$ & $(1,0)$\\
4 & $2/5$ & $z/16$ & $(0,0)$\\
5 & $1/4+z/64$ & $z/32$ & $(0,1)$\\
6 & $1/4+z/128$ & $1/12+z/64$ & $(0,1)$\\
7 & $13/48+z/256$ & $1/9+z/128$ & $(0,1)$\\
8 & $5/18+z/512$ & $55/432+z/256$ & $(0,1)$\\
9 & $487/1728+z/1024$ & $175/1296+z/512$ & $(1,1)$\\
10 & $2/5$ & $2161/15552+z/1024$ & \\
\bottomrule
\end{tabular}
\caption{The state after $t$ updates and the switches used in the next
update. }
\label{tab:itinerary}
\end{table}
\end{proof}
Figure~\ref{fig:orbit} shows the trajectory of the orbit from Theorem~\ref{thm:cycle}. 
\paragraph{Acknowledgment.} Generative AI models ChatGPT and Claude have been used for the refinement of the counterexample and \rev{the} exposition of the results.
\bibliographystyle{plain}
\bibliography{hk}

\end{document}